\documentstyle{amsppt} 
\pagewidth{12.5cm}\pageheight{19cm}\magnification\magstep1
\topmatter
\title{A study of intersections of Schubert varieties}\endtitle
\author M. Dyer and G. Lusztig\endauthor
\address{Department of Mathematics, University of Notre Dame, Notre Dame,
  IN 46556;
  Department of Mathematics, M.I.T., Cambridge, MA 02139}\endaddress
\endtopmatter   
\document
\define\lf{\lfloor}
\define\rf{\rfloor}

     \define\bco{\bar{\co}}

\define\bcb{\bar\cb}

\define\sqc{\sqcup}

\define\lb{\linebreak}

\define\part{\partial}

\define\m{\mapsto}
\define\do{\dots}

\define\lra{\leftrightarrow}

\define\T{\times}

\define\nl{\newline}
\redefine\i{^{-1}}

\define\ov{\overline}

\define\bbq{\bar{\QQ}_l}

\define\tr{\text{\rm tr}}

\redefine\d{\delta}

\define\io{\iota}

\define\p{\pi}

\define\CC{\bold C}

\define\QQ{\bold Q}

\define\ZZ{\bold Z}

\define\cb{\Cal B}

\define\ch{\Cal H}

\define\co{\Cal O}

\define\sha{\sharp}

\subhead 1\endsubhead
Let $G$ be a connected reductive group over $\CC$. Let $\cb$ the variety of Borel subgroups of $G$
and let $W$ be the Weyl group of $G$. Recall that the set of $G$-orbits
on $\cb\T\cb$ (with $G$ acting by simultaneous conjugation on the two
factors) is naturally in bijection $\co_w\lra w$ with $W$. It is well
known that for $w\in W$ the closure $\bco_w$ of $\co_w$ is
$\sqc_{w'\in W;w'\le w}\co_{w'}$; here $\le$ is the standard partial
order on $W$. Let $w\m|w|$ be the standard length function on $W$ and
let $w_0$ be the unique element of $W$ at which the length function
reaches its maximum. Let $(B^+,B^-)\in\co_{w_0}$. For $x\in W$ let

$\cb_x=\{B\in\cb;(B^+,B)\in\co_x\}$,
$\cb^x=\{B\in\cb;(B,B^-)\in\co_{x\i w_0}\}$.
\nl
It is well known that $\cb_x$ (resp. $\cb^x$) is isomorphic to
$\CC^{|x|}$ (resp. $\CC^{|x\i w_0|}$)
and that the closure of $\cb_x$ (resp. $\cb^x$) is

$\bcb_x:=\sqc_{x'\in W;x'\le x}\cb_{x'}=\{B\in\cb;(B^+,B)\in\bco_x\}$,
\nl
(resp. $\bcb^x:=\sqc_{x'\in W;x\le x'}\cb^{x'}
=\{B\in\cb;(B,B^-)\in\bco_{x\i w_0}\}$. )
\nl
For $x,y$ in $W$ we set

$\cb^y_x=\cb^y\cap\cb_x$.
\nl
This variety was introduced in \cite{KL79},
where it was shown that $\cb^y_x$ is nonempty if and only if $y\le x$.
Moreover, according to \cite{L98, 1.4}, if $y\le x$, then $\cb^y_x$ is
smooth of pure dimension $|x|-|y|$; according to \cite{R06, 7.1}, its
closure in $\cb$ is 
$$\bcb^y_x=\sqc_{(x',y')\in W\T W;y\le y'\le x'\le x}\cb^{y'}_{x'}.$$
Hence the intersection cohomology complex $K=IC(\bcb^y_x,\CC)$ is
defined. Here $\CC$ is viewed as a (constant) local system on $\cb^y_x$.

For $i\in\ZZ$ let $\ch^i(K)$ be the $i$-th cohomology sheaf of $K$. We
shall prove the following result.

\proclaim{Theorem 2} Assume that $y\le y'\le x'\le x$ and $i\in\ZZ$.
Then $\ch^i(K)|_{\cb^{y'}_{x'}}$ is a constant local system of rank say
$n^i_{y',x'}$. For $i$ odd we have $n^i_{y',x'}=0$. Moreover,
$$\sum_{j\in\ZZ}n^{2j}_{y',x'}q^j=P_{w_0y',w_0y}(q)P_{x',x}(q).\tag a$$
\endproclaim
In the case where $y=1$ we have $\bcb^y_x=\bcb_x$ and the theorem
specializes to the main result of \cite{KL80} which describes
the local intersection cohomology of $\bcb_x$ in terms of $P_{?,?}(q)$.
In fact, in no.3 we will deduce the theorem from this special case.
We have proved this theorem in 2003 but at that time did not write
down the proof. A proof was given in \cite{KWY13}. Since the proof in
\cite{KWY13}
seems to us more complicated than our proof, we thought that it might
be worth writing down our proof.

\subhead 3\endsubhead
We fix $y\le x$ in $W$.
We consider the diagram $\co_{w_0}@<\p'<<V'@>\io>>V@>\p>>\cb$ where 

$V=\{(B',B,B'')\in\cb\T\cb\T\cb;(B',B)\in\bco_x,(B,B'')\in\bco_{y\i w_0}\}$,

$V'=\{(B',B,B'')\in V;(B',B'')\in\co_{w_0}\}$,

$\p(B',B,B'')=B$, $\p'(B',B,B'')=(B',B'')$

and $\io$ is the obvious inclusion.

We have

$V=\sqc_{y',x'\text{ in }W;y\le y'\le x'\le x}V_{y',x'}$,

$V'=\sqc_{y',x'\text{ in }W;y\le y'\le x'\le x}V'_{y',x'}$

where
$$V_{y',x'}=\{(B',B,B'')\in\cb\T\cb\T\cb;(B',B)\in\co_{x'},
(B,B'')\in\co_{y'{}\i w_0}\},$$
$$\align&V'_{y',x'}=\{(B',B,B'')\in\cb\T\cb\T\cb;(B',B)\in\co_{x'},
(B,B'')\in\co_{w_0y'}, \\&
(B',B'')\in\co_{w_0}\}.\endalign$$

Note that $V_{y,x}$ (resp. $V'_{y,x}$) is a smooth open dense subset of
$V$ (resp. $V'$) hence the intersection cohomology complex
$K_V:=IC(V,\CC)$ (resp. $K_{V'}:=IC(V',\CC)$) is defined. (Here $\CC$ is
viewed as a (constant) local system on $V_{y,x}$ (resp. $V'_{y,x}$).
Assuming that $y\le y'\le x'\le x$ and $i\in\ZZ$ we show:

(a) {\it $\ch^i(K_V)|_{V_{y',x'}}$ is a constant local system of rank say
$m^i_{y',x'}$. For $i$ odd we have $m^i_{y'\le x'}=0$. Moreover,
$$\sum_{j\in\ZZ}m^{2j}_{y',x'}q^j=P_{w_0y',w_0y}(q)P_{x',x}(q).$$}

(b) {\it $\ch^i(K_{V'})|_{V'_{y',x'}}$ is a constant local system of
rank $m'{}^i_{y', x'}$. Moreover, we have
$m'{}^i_{y',x'}=m^i_{y',x'}$.}
\nl
Note that $V'$ is open in $V$ (via $\io$) and that
$V'_{y',x'}=\io\i(V_{y',x'})$. Hence (b) follows from (a).

Now $\p\i(B^+)$ can be identified with $\bcb_x\T\bcb_{w_0y}$
via $(B',B^+,B'')\m(B',B'')$.

(We use that $\bcb^y$ can be identified with $\bcb_{w_0y}$
via an automorphism of $G$ that exchanges $B^+,B^-$.)

Moreover, $\p'{}\i(B^+,B^-)$ can be identified
with $\bcb^y_x$ via $(B^+,B,B^-)\m B$. Under
these identifications $V_{y',x'}\cap\p\i(B^+)$ becomes the subset
$\cb_{x'}\T\cb_{w_0y'}$ of $\bcb_{x}\T\bcb_{w_0y}$
and $V'_{y',x'}\cap\p'{}\i(B',B'')$ becomes the subset $\cb^{y'}_{x'}$
of $\bcb^y_x$. 
Also, $\p$ and $\p'$ are locally trivial fibrations, compatible with
the $G$-actions (by conjugation on each factor) which are transitive on
their target. This implies that the local intersection cohomology of
$V$ has a simple relation to that of $\bcb_x\T\bcb_{w_0y}$
and that
the local intersection cohomology of $V'$ has a simple relation to that
of $\bcb^y_x$. In particular we see that (a) follows from the results
of \cite{KL80} applied to $\bcb_x$ and $\bcb_{w_0y}$ and that
$m'{}^i_{y',x'}=n^i_{y',x'}$, so that Theorem 2 follows from (a) and
(b).

\subhead 4\endsubhead
Let $q$ be an indeterminate.
For $i\in\ZZ$ we denote by $\bold H^i(\bcb^y_x)$
the $i$-th hypercohomology space of $\bcb^y_x$ with values in the complex
$K=IC(\bcb^y_x,\CC)$. The following result follows from Theorem 2
in the same way as Corollary 4.9 in \cite{KL80} followed from
\cite{KL80, 4.2, 4.3}. (We have proved this in 2003, unpublished; it was
also proved in \cite{P18} based on \cite{KWY13}.)

\proclaim{Corollary 5} For $i$ odd we have $\bold H^i(\bcb^y_x)=0$.
We have
$$\sum_{j\in\ZZ}\dim\bold H^{2j}(\bcb^y_x)q^j
=\sum_{y',x'\text{ in }W;y\le y'\le x'\le x}
P_{w_0y',w_0y}(q)R_{y',x'}(q)P_{x',x}(q)\tag a$$
where $P_{?,?}(q),R_{?,?}(q)$ are the polynomials in $\ZZ[q]$
defined in \cite{KL79}.
\endproclaim
Although the proof is standard, we will give it for completeness.
It is enough to prove the corollary when $G$ is replaced
by a connected reductive group (also denoted by $G$) with the same
Weyl group $W$ over 
 an algebraic closure of the finite field $F_p$ with $p$ elements
($p$ is a prime number) and the
 local system $\CC$ is replaced by $\bbq$ where $l$ is a prime
number $\ne p$. Let $F:G@>>>G$ be the Frobenius map corresponding to
a fixed split $F_p$-structure. Let $s\ge1$.
Note that $F^s$ induces a Frobenius map $\bcb^y_x@>>>\bcb^y_x$
(also denoted by $F^s$) 
and this induces automorphisms (denoted by $F^s$) of
each $\bold H^i(\bcb^y_x)$ and of the stalk $\ch^i(K)_z$ at any
$F^s$-fixed point $z$ of $\bcb^y_x$. By the Grothendieck-Lefschetz
fixed point formula we have
$$\align&\sum_i(-1)^i\tr(F^s,\bold H^i(\bcb^y_x))\\&=
\sum_{x',y' \text{ in }W;y\le y'\le x'\le x}
\sum_{z\in\cb^{y'}_{x'};F^s(z)=z}
\sum_{i\in\ZZ}(-1)^i\tr(F^s,\ch^i(K)_z).\endalign$$
From the proof of Theorem 2 and from
\cite{KL80, Thm.4.2}, we see that for $z\in\cb^{y'}_{x'}$ such that
$F^s(z)=z$, we have
$$\sum_{i\in\ZZ}(-1)^i\tr(F^s,\ch^i(K)_z)
=\sum_{j\in\ZZ}\tr(F^s,\ch^{2j}(K)_z)$$
and that $F^s$ acts on $\ch^{2j}(K)_z$ with only eigenvalues $p^{sj}$.
It follows that
$$\sum_{i\in\ZZ}(-1)^i\tr(F^s,\ch^i(K)_z)
=\sum_{j\in\ZZ}n_{y',x'}^{2j}p^{sj}.$$
Using this and Theorem 2 we see that
$$\align&\sum_i(-1)^i\tr(F^s,\bold H^i(\bcb^y_x))=\\&
\sum_{x',y'\text{ in } W;y\le y'\le x'\le x}
\sha\{z\in\cb^{y'}_{x'};F^s(z)=z\}P_{w_0y',w_0y}(p^s)P_{x',x}(p^s).
\endalign$$

We now use that for $x',y'$ as above we have
$\sha\{z\in\cb^{y'}_{x'};F^s(z)=z\}=R_{y',x'}(p^s)$
(a result of \cite{KL80}). We see that
$$\align&\sum_i(-1)^i\tr(F^s,\bold H^i(\bcb^y_x))=\\&
\sum_{x',y'\text{ in } W;y\le y'\le x'\le x}
P_{w_0y',w_0y}(p^s)R_{y',x'}(p^s)P_{x',x}(p^s).\tag b
\endalign$$

In particular, the left hand side of (b) is a
polynomial in $p^s$ when $s$ varies in $\{1,2,3,\do\}$.
By Deligne's  theorem, for any $i\in\ZZ$, any eigenvalue
of $F^s$ on $\bold H^i(\bcb^y_x)$ has absolute value $p^{is/2}$ after
applying to it any isomorphism of $\bbq$ with $\CC$. Using this and (b)
we see that
$\bold H^i(\bcb^y_x)=0$ when $i$ is odd and that, when $j\in\ZZ$,
any eigenvalue of $F^s$ on $\bold H^{2j}(\bcb^y_x)$ is equal to $p^{js}$.
Thus, the left side of (b) is equal to
$$\sum_{j\in2\ZZ}\dim\bold H^{2j}(\bcb^y_x)p^{js},$$
so that the corollary follows from (b).

\subhead 6\endsubhead
In the rest of this paper, $W$ is any Coxeter group with standard length function
$w\m|w|$ and standard partial order $\le$. For $y\le x$ in $W$
the polynomials
$P_{y,x},R_{y,x}$ in $q$ are defined as in \cite{KL79}; moreover the
polynomial $Q_{y,x}$ in $q$ is defined as in \cite{KL80,\S2}. According
to \cite{KL79}, when $W$ is a Weyl group we have $Q_{y,x}=P_{w_0x,w_0y}$.
Hence in this case the right hand side of 5(a) can be written as
$$\Xi_{y,x}(q)=\sum_{x',y'\text{ in }W;y\le y'\le x'\le x}
Q_{y,y'}(q)R_{y',x'}(q)P_{x',x}(q).\tag a$$
Note that in this form $\Xi_{y,x}(q)\in\ZZ[q]$ is well defined
for any $y\le x$ in any Coxeter group $W$. 

Let $\bar{}:\ZZ[q,q\i]@>>>\ZZ[q,q\i]$ be the ring involution which takes
$q$ to $q\i$ and $q\i$ to $q$. Recall from \cite{KL80} that
$$\ov{P_{x',x}(q)}
=\sum_{x'\le u\le x}R_{x',u}(q)P_{u,x}(q)q^{-|x|+|x'|},$$
$$\ov{R_{y',x'}(q)}=R_{y',x'}(q)q^{-|x'|+|y'|}(-1)^{|x'|+|y'|},$$
$$\ov{Q_{y,y'}(q)}
=\sum_{y\le v\le y'}Q_{y,v}(q)R_{v,y'}(q)q^{-|y'|+|y|}.$$
Using these in (a) gives
$$\Xi_{y,x}(q)=q^{|x|}\sum_{y'\text{ in }W;y\le y'\le x}
q^{-|y'|}Q_{y,y'}(q)\ov{P_{y',x}(q)}, \tag b$$
$$\Xi_{y,x}(q)=q^{-|y|}\sum_{x'\text{ in }W;y\le x'\le x}
q^{|x'|}\ov{Q_{y,x'}(q)}P_{x',x}(q).\tag c$$
Replacing $y'$ by $x'$ in the formula (b) and comparing with (c) shows
that
$$\ov{\Xi_{y,x}(q)}=q^{|y|-|x|}\Xi_{y,x}(q).\tag d$$
(If $W$ is a Weyl group this follows from 5(a) using Poincar\'e
duality.)

In (a), standard degree bounds for $Q_{y,y'}(q)$, $R_{y',x'}(q)$ and  $P_{x',x}(q)$ imply  that $\Xi_{y,x}$ has the same degree and leading coefficient as $R_{y,x}$. Therefore

(e) {\it $\Xi_{y,x}$ is a  monic polynomial in $q$ of  degree $|x|-|y|$.}

It follows immediately from (b)--(c), using \cite{KL80, 2.1.6}, that 
$$\sum_{z\in W;y\le z\le x}(-1)^{|y|+|z|}\Xi_{y,z}(q)\Xi_{z,x}(q)=
\d_{x,y}.\tag f$$

\subhead 7 \endsubhead
Recall from \cite{KL79} the definition of the generic Iwahori-Hecke algebra over $\ZZ[q^{1/2},q^{-1/2}]$, with generators $T_{r}$ for  simple reflections $r$ subject to the braid relations of $W$ and quadratic relations $(T_{r}+1)(T_{r}-q)=0$.
It has a standard basis $(T_{w})_{w\in W}$ and two bases
$(C_{w})_{w\in W}$ and $(C'_{w})_{w\in W}$ defined in \cite{KL79}. From \cite{KL79,(1.1.b),(1.1.c)} and
\cite{KL80, 2.1.6}, one directly calculates using 6(c) that
$$C_{x}'=\sum_{y\in W;y\leq x}q^{(|y|-|x|)/2}\Xi_{y,x}C_{y}.
\tag a$$
Alternative proofs of 6(d), 6(e) and 6(b) may be based on this formula. 
 
 \subhead 8 \endsubhead
 By 6(d), for $y\leq x$ in $W$, we may write  $\Xi_{y,x}=\sum_{i=0}^{N}a_{i}q^{i}$ where  each $a_{i}\in \ZZ$ and  $N=|x|-|y|$. We have  $a_{i}=a_{N-i}$ for $i=0,\ldots, N$ by 6(d) and $a_{0}=a_{N}=1$ by 6(e).
 By \cite{EW14}, the polynomials $P_{?,?}$
 have non-negative coefficients. We expect that the polynomials
 $Q_{?,?}$ also  have non-negative coefficients.
 This would imply by 6(b) or 6(c) that  one has $a_{i}\geq 0$ for all
 $i=0,\ldots, N$.
 
We expect in general that
$$0\leq a_{0}\leq a_{1}\leq\ldots \leq a_{\lf \frac{N+1}{2}\rf}. \tag{a}$$ When $W$ is a Weyl group, this follows from Corollary 5 using the
hard Lefschetz theorem in intersection cohomology. 

More generally, if $W$ is crystallographic, then
$\cb^y_x=\cb^y\cap\cb_x$ can still
be defined as in \cite{KL80,5.2} (where $W$ is assumed to be an affine
Weyl group, but this assumption is unnecessary) and, as in
{\it loc.cit.}, one
can consider the projective variety $\bcb^y\cap\bcb_x$. It is likely
that the analogues of Theorem 2 and Corollary 5
(with $P_{w_0y',w_0y}$ replaced by $Q_{y,y'}$)
hold for this variety (with a similar
proof), so that the positivity statement above would hold in this case.

When $W$ is a finite Coxeter group, (a)  follows  from 
7(a), (the proof of) \cite{DL90, (2.7)(ii)} and \cite{EW21, (3)}.

 \subhead 9 \endsubhead
 We wish to thank Thomas Lam  for providing to us the reference
 to \cite{P18} (after we posted a first version of this paper)
 which then led us to \cite{KWY13}.

\enddocument